\documentclass[11pt]{article}
\usepackage{amsfonts, amssymb, mathtools}
\usepackage[a4paper,top=5.5cm,bottom=4cm,left=2.5cm, right=2.5cm,foot=1cm]{geometry}
\usepackage{authblk}
\usepackage{bm}
\ifx\pdftexversion\undefined
 \usepackage[dvips]{graphicx}
\else
  \usepackage[pdftex]{graphicx}
\fi
\usepackage{hyperref}
\usepackage{tikz}
\usetikzlibrary {automata}
\newcommand{\emailaddress}[1]{{\sf#1}}
 \let\LaTeXtitle\title
\renewcommand{\title}[1]{\LaTeXtitle{\Large{\textbf{#1}}}}
\usepackage{xcolor}

\title{Reconstruction of Transient Anisotropic Diffusion in Transient Diffusion Models} 

\date{\vspace{-6ex}} 

\author[1]{\underline{Luchesi, Vanda M.}} 
\affil[1]{Federal University of Santa Catarina, \emailaddress{v.luchesi@ufsc.br}} 

\DeclareMathOperator*{\argmin}{argmin}
\newcommand{\R}{\mathbb{R}}

\usepackage{float}
\usepackage{multirow}
\usepackage{multicol}
\usepackage{subfig}
\usepackage{booktabs}
\usepackage{authblk}
\usepackage{tikz}

\tikzstyle{vertex} = [draw,circle,fill=white,thick,inner sep=0pt,minimum size=5pt]
\tikzstyle{dots} = [draw,circle,fill=black,inner sep=0pt,minimum size=1pt]

\begin{document}

\maketitle


This presentation abstract, presented at 12th Applied Inverse Problems Conference, to be held at FGV EMAp, Rio de Janeiro, Brazil, in 2025, describe the theory results of applying the Chebyshev pseudospectral method (CPM) to reconstruct the anisotropic transient diffusion in a transiente two-dimensional diffusion model with boundary conditions. This model can be easily applied to a wide range of diffusion phenomena across various science domains.


The diffusion problem within a finite domain  $(0,t_f] \times \Omega$,  $t_f >0$,  with $\Omega=(0,1)\times (0,1) \subseteq  \mathbb{R}^2$ can be apresented by equation
\begin{equation}\label{eq:heat}
U_t (t,x,y)=\nabla \cdot [{K}(t,x,y)\nabla 
U(t,x,y)] + g(t,x,y)
\end{equation}
in the space-time  $(0,t_f] \times \Omega$, where the time-space-dependent tensor diffusion given as
\begin{equation}\label{eq:diffusion}
K(t,x,y) = \left[ \begin{array}{cc}
k_{11} (t,x,y) & k_{12}(t,x,y) \\ 
k_{21}(t,x,y) & k_{22} (t,x,y)
\end{array}  \right]
\end{equation}
with $g(t,x,y)$ a time dependent source. 
$K(t,x,y)$ for $t\neq 0$ regarded as unknown and has to be estimated from temperature data or the pressure in a transport problem in porous media. Many studies to identify tensor diffusion (or conductivity) have been developed and reported in the literature 
\cite{Bo_Luchesi_Bazan} \cite{B_Luchesi_Bazan} and \cite{Bahssini}. 
In this sense, we will impose that $K(t,x,y)$ will be positive definite matrix with $k_{11} (t,x,y) k_{22}(t,x,y)-k_{21}(t,x,y) k_{12} (t,x,y) >0$ ensured the operations are well-defined and stable.
 \begin{figure}[!ht]
 \centering
\includegraphics[scale=0.3]{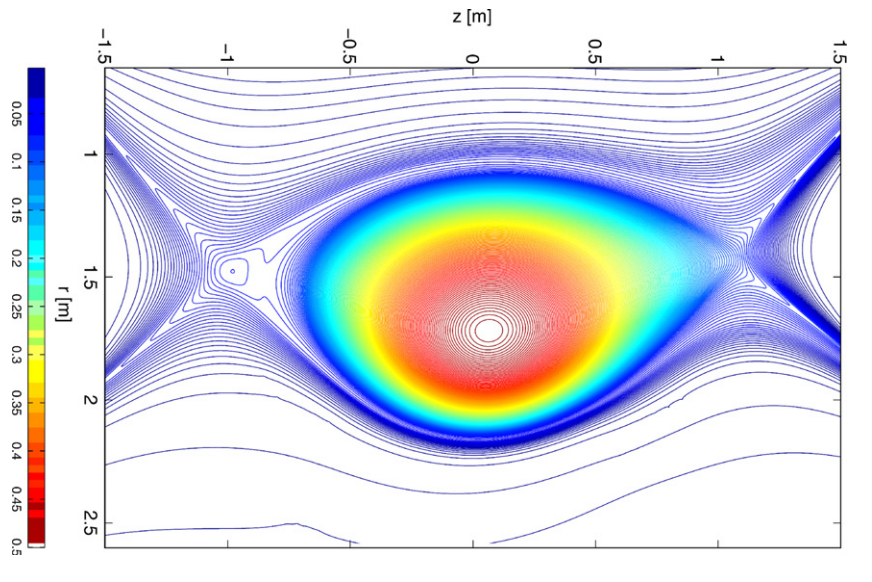}
 \includegraphics[scale=0.4]{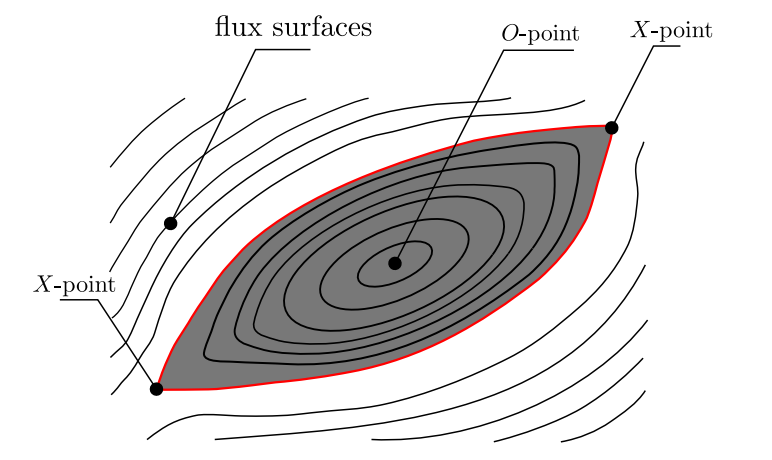}
 \caption{Magnetic island in Plasma,  and Scheme of flux surface for identification with elliptical coordenates (sources: \cite{Gunter2007} and \cite{Thesis} respectively)}
  \end{figure}

for a two-dimensional problem, the diffusion tensor by $K(t,x,y) =\Phi(t)\varLambda(t) \Phi(t) $, with the time dependent rotation matrix $\Phi(t)$ given as
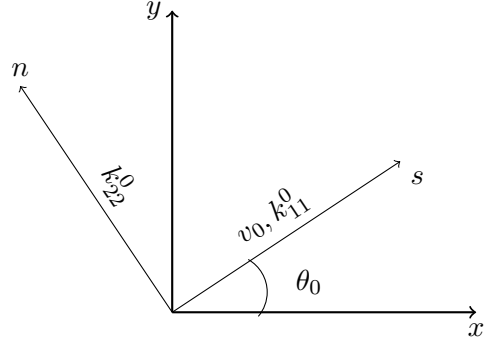
\begin{figure}[!ht]
	\centering

\begin{tikzpicture}
    \draw[->, thick] (0,0) -- (4,0) node[below] {$x$};
    \draw[->, thick] (0,0) -- (0,4) node[left] {$y$};
    
    \draw[->] (0,0) -- (-2,3) node[sloped,midway,above] {$k^0_{22}$} node[above] {$n$};
    \draw[->] (0,0) -- (3,2) node[sloped,midway,above] {$v_0,k^0_{11}$} node[below right] {$s$};
    
    \draw (1,0.7) arc[start angle=60,end angle=-40,radius=0.5cm];
    \node at (1.8,0.4) {$\theta_0$}; 
  \end{tikzpicture}
 \caption{ Ilustration of the change of variable given by director cosseno of isothermal plane: $k^0_{11}=k^0_{11}(x_0, y_0)$ e $k^0_{22}=k^0_{22}(x_0, y_0)$,  are the principal diffusion for a fixed $t_0 \in [0,t_f]$ (paralell and ortogonal respectively) with  $v_0=(v_1(t_0),v_2(t_0))^T=(\cos \theta(t_0), \sin \theta(t_0))^T$ and $\theta(t_0)=\theta_0$}
 \label{change}
\end{figure}
$$\Phi(t) =
\begin{bmatrix}
\cos \theta (t) & -\sin \theta (t) \\
\sin \theta (t ) & \cos \theta (t)
\end{bmatrix}$$
where $\theta (t)$ is the tetha function which represents the angular function varing with the time, see fig. 

$$\varLambda(x,y)=\left[ \begin{array}{cc}
k_{11}(x,y)  & 0 \\ 
0 & k_{22}(x,y) 
\end{array}  \right]$$
is the diagonal matrix where $k_{ii}(x,y)$ are the principal diffusion (paralell and ortogonal) given by director cosseno of isothermal plane.
Here, we propose the following new approach, time dependent, for anisotropic diffusion:
\begin{equation}\label{eq:diffusion2}
K(t,x,y) =\left[ \begin{array}{cc}
k_{11}(x,y)\cos^2( \theta (t))+k_{22}(x,y)\sin^2(\theta(t)) & (k_{11}(x,y)-k_{22}(x,y))\sin (\theta (t))\cos (\theta (t)) \\
(k_{11}(x,y)-k_{22}(x,y))\sin (\theta (t))\cos (\theta (t)) & k_{11}(x,y)\sin^2( \theta (t))+k_{22}(x,y)\cos^2(\theta(t))
\end{array} \right],
\end{equation} 
Considering $\sin (\theta (t))\cos (\theta (t))= \frac{1}{2}\sin (2\theta (t))$
\begin{equation}\label{eq:diffusion3}
K(t,x,y) =\left[ \begin{array}{cc}
k_{11}(x,y)\cos^2( \theta (t))+k_{22}(x,y)\sin^2(\theta(t)) & \frac{1}{2}(k_{11}(x,y)-k_{22}(x,y))\sin (2\theta (t)) \\
\frac{1}{2}(k_{11}(x,y)-k_{22}(x,y))\sin (2\theta (t)) & k_{11}(x,y))\sin^2( \theta (t))+k_{22}(x,y)\cos^2(\theta(t))
\end{array} \right],
\end{equation}

This formulation for  $\theta(t)=0$, began:
\begin{itemize}
 \item Ortrotopic diffusion,  
 \begin{equation}\label{eq:dif}
K(0,x,y)  =\left[ \begin{array}{cc}
k_{11}(x,y) & 0 \\
0& k_{22}(x,y)
\end{array} \right]
\end{equation} 
 In this context, this formulation of anisotropy can be seeing also by related to diffusion coefficients, where the parallel diffusion coefficient is significantly greater than the perpendicular diffusion coefficient $(k_{11}(x,y) \gg k_{22}(x,y)$).
 \item Isotropic diffusion when $k_{11}(x,y)=k_{22}(x,y)=k(x,y)$.
 \begin{equation}\label{eq:dif2}
K(0,x,y) =\left[ \begin{array}{cc}
k(x,y)  & 0 \\
0 & k(x,y)
\end{array} \right],
\end{equation}
Here, this formulation concerning also with temperature distribution, in which the perpendicular temperature gradient far exceeds the parallel temperature gradient $(\nabla_{22} U \gg \nabla_{11}U$).
\end{itemize}

\section{Anisotropic diffusion}
We will rewrite the eq. \ref{eq:heat} as
\begin{equation}\label{eq:heat2}
U_t =\frac{\partial }{\partial x}\left[ k_{11}(t,x,y) U_x+k_{12}(t,x,y) U_y\right] +\frac{\partial }{\partial y} \left[ k_{21}(t,x,y)) U_x+k_{22}(t,x,y) U_y\right] 
 + g(t,x,y)
\end{equation}
Replacing the $k_{ii}(x,y,t)$ by given in eq. \ref{eq:diffusion2} we can write \ref{eq:heat2} as
\begin{eqnarray}\label{eq:heat4}
U_t = 
\frac{\partial }{\partial x}\left( k_{11}(x,y)\cos^2( \theta (t))+k_{22}(x,y)\sin^2(\theta(t)) U_x+(k_{11}(x,y)-k_{22}(x,y))\sin (\theta (t))\cos (\theta (t)) U_y\right) + \nonumber \\ \frac{\partial }{\partial y}\left( (k_{11}(x,y)-k_{22}(x,y))\sin (\theta (t))\cos (\theta (t) U_x+k_{11}(x,y)\sin^2( \theta (t))+k_{22}(x,y)\cos^2(\theta(t))U_y\right) + \nonumber \\
 + g(t,x,y)
\end{eqnarray}
Considering $\frac{1}{2}\sin (2\theta (t))= \sin (\theta (t))\cos (\theta (t)$ the eq. \ref{eq:heat4} return the following equation
\begin{eqnarray}\label{eq:heatF}
 U_t = 
\cos^2( \theta (t))\frac{\partial }{\partial x}\left( k_{11}(x,y)U_x\right)+\sin^2(\theta(t)) \frac{\partial }{\partial x}\left( k_{22}(x,y)U_x\right)+\nonumber \\ +\frac{1}{2}\sin (2\theta (t))\frac{\partial }{\partial x}\left(k_{11}(x,y) U_y \right)- \frac{1}{2}\sin (2\theta (t))\frac{\partial }{\partial x}\left(k_{22}(x,y) U_y \right)  + \nonumber \\+\frac{1}{2}\sin (2\theta (t))\frac{\partial }{\partial y}\left(k_{11}(x,y) U_x \right)-\frac{1}{2}\sin (2\theta (t))\frac{\partial }{\partial y}\left(k_{22}(x,y) U_x \right)+ \nonumber \\ +\sin^2( \theta (t))\frac{\partial }{\partial y}\left(k_{11}(x,y)U_y \right)+\cos^2( \theta (t))\frac{\partial }{\partial y}\left(k_{22}(x,y)U_y \right)  + g(t,x,y)
\end{eqnarray}
The inicial conditions is.
\begin{align}
\label{eqe:Initial}
U(0,x,y)=U_0 (x,y) &\quad \text{in} \ [0,1] \times [0,1].
\end{align}
But, for $(t,x,y) \in [0,1] \times [0,1] \times [0,t_f]$, along with boundary conditions for all $t\in (0,t_f]$ we will make the following change. 

The first boundary condition given by
\begin{align}
\label{eqe:boundary1}
k_{11}(t,0,y)U_x(t,0,y)+k_{12}(t,0,y)U_y(t,0,y)=U(t,0,y)-f_1(t,0,y), &\quad  y \in (0,1),
\end{align}
with  the approach $$k_{11}(t,0,y)=k_{11}(0,y)\cos^2( \theta (t))+k_{22}(0,y)\sin^2(\theta(t))$$
$$k_{12}(t,0,y)=\frac{1}{2}(k_{11}(0,y)-k_{22}(0,y))\sin (2\theta (t)) $$
Become, for $0< \theta(t) < \pi/2$
\begin{align}\notag
\cos^2 \theta(t)k_{11}(0,y)  U_x(0,y,t)= 
U(t,0,y)-f_1(t,0,y)-k_{22} (0,y)\sin^2 (\theta (t))U_x(t,0,y)\\ -\frac{1}{2}[k_{11}(0,y)-k_{22}(0,y)] 
\sin(2 \theta(t)) U_y(t,0,y)  
\end{align} or, when necessary
\begin{align}\notag
\sin^2 \theta(t)k_{22}(0,y)  U_x(t,0,y)= 
U(t,0,y)-f_1(t,0,y)-k_{11} (0,y)\cos^2 (\theta (t))U_x(t,0,y)\\ -\frac{1}{2}[k_{11}(0,y)-k_{22}(0,y)] 
\sin(2 \theta(t)) U_y(t,0,y)  
\end{align} 
 The second boundary conditions 
\begin{align}
\label{eqe:boundary2}
k_{11}(t,1,y)U_x(t,1,y)+k_{12}(t,1,y)U_x(t,1,y)=- U(t,1,y) )+f_2(t,1,y),  &\quad  \ y \in (0,1),
\end{align}
with the approach
$$k_{11}(t,1,y)=k_{11}(1,y)\cos^2( \theta (t))+k_{22}(1,y)\sin^2(\theta(t))$$
$$k_{12}(t,1,y)=\frac{1}{2}(k_{11}(1,y)-k_{22}(1,y))\sin (\theta (t))\cos (\theta (t))$$
Become for $0< \theta(t) < \pi/2$
\begin{align}\notag
k_{11}(1,y)  U_x(1,y,t)= 
\frac{1}{\cos^2 \theta(t)} \left \{- U(t,1,y)+f_2(t,1,y)-\frac{1}{2}[k_{11}(1,y_j)-k_{22}(1,y)] \right. \\ \left. \sin(2 \theta(t)) U_y(1,y,t)-k_{22} (1,y)\sin^2 (\theta (t))U_x(t,1,y) \right\}
\end{align}
or, when necessary
\begin{align}\notag
\sin^2 \theta(t)k_{22}(1,y)  U_x(t,1,y)= 
U(t,1,y)-f_2(t,1,y)-k_{11} (1,y)\cos^2 (\theta (t))U_x(t,1,y)\\ -\frac{1}{2}[k_{11}(1,y)-k_{22}(1,y)] 
\sin(2 \theta(t)) U_y(t,0,y)  
\end{align} 
The third boundary conditions, given by 
\begin{align}
\label{eqe:boundary3}
k_{21}(t,x,0)U_x(t,x,0)+k_{22}(t,x,0)U_y(t,x,0)=U(t,x,0)-f_3(t,x,0),   &\quad  x \in (0,1),
\end{align}
with the approximation
$$k_{21}(t,x,0)=\frac{1}{2}(k_{11}(x,1)-k_{22}(x,1))\sin (2\theta (t))$$
$$k_{22}(t,x,0)=k_{11}(x,y))\sin^2( \theta (t))+k_{22}(x,y)\cos^2(\theta(t))$$
Become, for $0< \theta(t) < \pi/2$
\begin{align}\notag
k_{22}(x,0)  U_y(t,x,0)= 
\frac{1}{\cos^2 \theta(t)}  \left\{ U(t,x,0)-f_3(t,x,0)-\frac{1}{2}[k_{11}(x,0)-k_{22}(x,0)] \times \right. \\ \left. \sin(2 \theta(t)) U_x(t,x,0)  -k_{11} (x,0)\sin^2 (\theta (t))U_y(t,x,0) \right\}, & \quad  x \in (0,1),
\end{align}
or, when necessary
\begin{align}\notag
k_{11}(x,0)  U_y(t,x,0)= 
\frac{1}{\sin^2 \theta(t)}  \left\{ U(t,x,0)-f_3(t,x,0)-\frac{1}{2}[k_{11}(x,0)-k_{22}(x,0)] \times \right. \\ \left. \sin(2 \theta(t)) U_x(t,x,0)  -k_{22} (x,0)\cos^2 (\theta t)U_y(t,x,0) \right\}, & \quad  x \in (0,1),
\end{align}
The fourth boundary conditions, given by 
\begin{align}
\label{eqe:boundary4}
k_{21}(t,x,1)U_x(t,x,1)+k_{22}(t,x,1)U_y(t,x,1) =- U(t,x,1)+f_4(t,x,1),   &\quad  x \in (0,1),
\end{align}

with the approximation
$$k_{21}(t,x,1)=\frac{1}{2}(k_{11}(x,1)-k_{22}(x,1))\sin (2\theta (t))$$
$$k_{22}(t,x,1)=k_{11}(x,1))\sin^2( \theta (t))+k_{22}(x,1)\cos^2(\theta(t))$$
Become
\begin{align}\notag
k_{22}(x,1)  U_y(t,x,1)= 
\frac{1}{\cos^2 \theta(t)} \left \{- U(t,x,1)+f_4(t,x,1)-\frac{1}{2}[k_{11}(x,1)-k_{22}(x,1)] \times \right. \\ \left. \sin(2 \theta(t)) U_x(t,x,1)-k_{11}(x,1))\sin^2( \theta (t))U_y(t,x,1) \right\}&\quad  x \in (0,1),
\end{align}
or, when necessary
\begin{align}\notag
k_{11}(x,1)  U_y(t,x,1)= 
\frac{1}{\sin^2 \theta(t)} \left \{- U(t,x,1)+f_4(t,x,1)-\frac{1}{2}[k_{11}(x,1)-k_{22}(x,1)] \times \right. \\ \left. \sin(2 \theta(t)) U_y(t,x,1)-k_{22}(x,1))\cos^2( \theta (t))U_x(t,x,1) \right\}&\quad  x \in (0,1),
\end{align}
\section{Chebyshev discretization for spatial variable}
To solve this problem, the PDE in Eq. \ref{eq:heat} will be transformed into a time dependent system of ordinary differential equations (ODEs) applying the CPM method. And later, the ODEs system will be solved by other numerical method. Firstly, the mesh consisting of $(n+1)\times(n+1)$  points on the domain will be transformed on $(n+1)$ Chebyshev Gauss-Lobatto $x_i$ and $y_i$ in the  horizontal and vertical directions, respectively. 
Then spatial derivatives on a mesh points will be aproximatted on Chebyshev Gauss-Lobatto points defined by
\begin{equation}
 \hat{x}_{i} = \frac{1}{2}\left[1-cos(i\pi/n)\right], \quad i= 0, 1, n
\end{equation}
\begin{equation}
 \hat{y}_{j} = \frac{1}{2}\left[1-cos(j\pi/n)\right], \quad j= 0, 1, n
\end{equation}
with the approximation being based on a matrix-vector product of the Chebyshev differentiation matrix with a vector of function values on the grid.
Let $D$  the $(n+1)\times (n+1)$ Chebyshev differentiation matrix and assume that it is decomposed as
\begin{equation}\label{twoform0}
D = \left  [\begin{array}{c}
r_0^T\\
r_1^T\\
\vdots\\ r_n^T 
\end{array} \right ] =  \left [ d_0,d_1,\dots, d_n\right ], \quad r_i, d_i \in \mathbb{R}^{n+1}, \; i=0,\dots, n.
\end{equation}
The spatial derivatives in (\ref{eq:heatF}) will be transform into a discrete form using the properties of $D$, for this we will first considering the derivative with respect to $x$. Identifying the vectors with the following notations
\begin{equation}\label{Tj}
{\bf U}_j (t) = [U(t,x_0,y_j), U(t,x_1,y_j), \dots, U(t,x_n,y_j)]^T, \;\; j=0,1,\dots, n,
\end{equation}
For $z=x$ and  $z=y$
\begin{equation}\label{T-derj}
{\bf U}^z_j (t)=\left  [ \dfrac{\partial U(t,x_0,y_j)}{\partial z },  \dfrac{\partial U(t,x_1,y_j)}{\partial z}, \dots,
\dfrac{\partial U(t,x_n,y_j)}{\partial z}\right ]^T, \;\; j=0,1,\dots, n.
\end{equation}
and next identification show that there are no constraints after CPM method applied

\begin{equation}\label{whole}
[k_{11}^j {\bf U}_j^{x}]^x :=
\left [ \begin{array}{c} 
\dfrac{\partial}{\partial x}\left ( k_{11}(x_0,y_j)
\dfrac{\partial U(t,x_0,y_j)}{\partial x}\right ) \\
\vdots \\
\dfrac{\partial}{\partial x}\left (  k_{11} (x_n,y_j)
\dfrac{\partial U(t,x_n,y_j)}{\partial x}\right )  \end{array}\right ]\approx
D ({\bf K}^1_{j}{{{\bf U}}^x_j}),\;\ j=0,\dots,n,
\end{equation}
where ${\bf K}^1_{j}$  denotes  the matrix diffusion given by
\begin{equation}\label{Kj}
{\bf K}^1_{j} = {\rm diag}(k_{11}(x_0,y_j), k_{11}(x_1,y_j) , \dots, k_{11}(x_n,y_j)),\;\; j=0,1,\dots, n.
\end{equation}

and for $j=1,\dots, n-1$ we have
$${{\bf U}}^x_j \approx D_2{\bf U}_j$$ 
\subsection{Boundary Condition}

Applying now, the CPM method, for the boundary conditions
Eqs.((\ref{eqe:boundary1})--(\ref{eqe:boundary2}),  
at $(x_0,y_j)$, $j=0,\dots, n$, and
at $(x_n,y_j)$, $j=0,\dots, n$. 
\begin{enumerate}
 \item For arbitrary $j$, using such boundary conditions
and similar aproximatting given in Eq.(\ref{whole})  we have
\begin{align}\notag
[k_{11}^j {\bf U}_j^{x}]^x &\approx
\sum_{j=0}^n d_jk_{11}(x_i,y_j) U_x(t,x_i,y_j) \\\notag
&= d_0k_{11}(0,y_j) U_x(0,y_j,t) + 
\displaystyle \sum_{j=1}^{n-1} d_jk_{11}(x_i,y_j) 
 U_x(t,x_i,y_j) 
 +d_nk_{11}(1,y_j)  U_x(t,1,y_j)\\ 
&= \left [ d_0 \beta_j^0e_1^T +D_1\check{\bf K}_j^{1}D_2 -d_n \beta_j^ne_{n+1}^T\right ]^T
{\bf U}_j + {\bf b}_j\\ \label{bott-top}
&\doteq {\bf A}_j{\bf U}_j + {\bf b}_j,
\end{align}
$${\bf b}_j =\frac{1}{\cos^2 \theta(t)} \{ -d_0f_1(t,0,y_j)+ d_nf_2(t,1,y_j)\},$$
\begin{align}\notag
{ \beta}^0_j= 
\frac{1}{\cos^2 \theta(t)}  \left\{1-\frac{1}{2}[k_{11}(0,y_j)-k_{22}(0,y_j)] \sin(2 \theta(t))  -k_{22} (0,y_j)\sin^2 (\theta (t) \right\}
\end{align}
\begin{align}\notag
{\bf \beta}^n_j= 
\frac{1}{\cos^2 \theta(t)}  \left\{1 -\frac{1}{2}[k_{11}(1,y_j)-k_{22}(1,y_j)]  \sin(2 \theta(t)) -k_{22} (1,y_j)\sin^2 (\theta (t)) \right\}
\end{align}

\begin{align}
{\bf A}^1_j &= \left [d_0 {\beta}^0_je_1^T +D_1\check{\bf K}_j^{1}D_2 -d_n {\beta}^n_j e_{n+1}^T\right ]^T 
\end{align}
with
\begin{equation*}
D_1 = [d_1,\dots,d_{n-1}], \quad  D_2 = \left [ \begin{array}{c}
r_1^T\\
\vdots \\
r_{n-1}^T\\
\end{array} \right ], 
\end{equation*}

\item Similary, now for arbitrary $i$, applying CPM for eqs.(\ref{eqe:boundary3})--(\ref{eqe:boundary4}),  
at $(x_i,y_0)$, $i=0,\dots, n$, and
at $(x_i,y_n)$, $i=0,\dots, n$, we obtain
$$[k_{22}^i {\bf U}_i^{y}]^y\doteq {\bf A}^2_i{\bf U}_i +{\bf b}^y_i$$
with
\begin{align}
{\bf A}^2_i &= d_0{\beta}^0_ie_1^T +D_1\check{\bf K}^2_{i}D_2 +d_n {\beta}^n_i e_{n+1}^T \quad \text{and}\\
{\bf b}^y_i &= \frac{1}{\cos^2 \theta(t)} \{-d_0f_3(t,x_i,0)+d_n f_4(t,x_i,1)\},
\end{align}
\begin{align}
{\bf \beta}^0_i=\frac{1}{\cos^2 \theta(t)}  \left\{ 1-\frac{1}{2}[k_{11}(x_i,0)-k_{22}(x_i,0)]  \sin(2 \theta(t))   -k_{11} (x_i,0)\cos^2 (\theta t)
\right\}
\end{align}
\begin{align}
{\bf \beta}^n_i=\frac{1}{\cos^2 \theta(t)} \left \{1-\frac{1}{2}[k_{11}(x,1)-k_{22}(x_i,1)]  \sin(2 \theta(t)) -k_{11}(x_i,1))\sin^2( \theta (t)) \right\}
\end{align}
\begin{equation}\label{K2im}
 \check{\bf K}^2_{i} = {\rm diag}( k_{22}(x_i,y_1), \dots, k_{22}(x_i,y_{n-1}),\;\; i=0,1,\dots, n.
\end{equation}
Because for $i=1,\dots, n-1$ we have
$${{\bf U}}^y_i \approx D_2{\bf U}_i$$ 

\item Similary, written 
\begin{align}\notag
[k_{22}^i {\bf U}_i^{y}]^x &\approx
\sum_{i=0}^n d_ik_{22}(x_i,y_j)  \dfrac{\partial 
	U(x_i,y_j,t)}{\partial y } \\\notag
&= d_0k_{22}(x_i,0)  \dfrac{\partial U(x_i,1,t)}{\partial y } + 
\displaystyle \sum_{i=1}^{n-1} d_ik_{22}(x_i,y_j) 
\dfrac{\partial U(x_i,y_j,t)}{\partial y } 
 +d_nk_{22}(x_i,1)  \dfrac{\partial U(x_i,1,t)}{\partial y } \\\notag 
 &= \left [ d_0 \beta^0_ie_1^T + D_1\check{\bf K}^2_iD_2 -d_n\beta^n_ie_{n+1}^T\right ] 
{\bf U}_i+\frac{1}{\cos^2 \theta(t)} \{-d_0f_3(x_i,y_j,t) 
 + d_nf_4(x_i,y_j,t) \}\\ \label{bott-top}
\end{align}
we can get
\begin{equation}\label{finite-yx}
[k_{22}^i {\bf U}_i^{y}]^x=\dfrac{\partial}{\partial x}\left ( k_{22}(x_i,y_j)\dfrac{\partial U(x_i,y_j,t)}{\partial y}\right ) \approx
D ({\bf K}^2_{i}{{{\bf U}}^y_i})
\doteq  {\bf A}^{2}_i{\bf U}_i + {\bf b}^{y}_i, \, i=0,\dots, n.
\end{equation}

\begin{align}
{\bf A}^{2}_i & = d_0\beta^0_ie_1^T +D_1 \check{\bf K}^2_{i}D_2 -d_n\beta^n_ie_{n+1}^T \, \,  \quad \text{and}\\
{\bf b}^y_i &= \frac{1}{\cos^2 \theta(t)} \{-d_0f_3(t,x_i,0)+ d_nf_4(t,x_i,1)\},
\end{align}
\item Furthemore, write
\begin{align}\notag
[k_{11}^j {\bf U}_j^{x}]^y &\approx
\sum_{j=0}^n d_jk_{11}(x_i,y_j)  \dfrac{\partial 
U(x_i,y_j,t)}{\partial x } \\\notag
&= d_0k_{11}(x_0,y_j)  \dfrac{\partial U(x_0,y_j,t)}{\partial x } + 
\displaystyle \sum_{j=1}^{n-1} d_jk_{11}(x_i,y_j) 
\dfrac{\partial U(x_i,y_j,t)}{\partial x } 
 +d_nk_{11}(x_n,y_j)  \dfrac{\partial U(x_n,y_j,t)}{\partial x} \\\notag 
 &= \left [ d_0\beta^0_j e_1^T + D_1\check{\bf K}^1_jD_2 -d_n\beta^n_je_{n+1}^T\right ] 
{\bf U}_j+\frac{1}{\cos^2 \theta(t)} \{ -d_0f_1(x_i,y_j,t) 
 + d_nf_2(x_i,y_j,t \}\\ \label{bott-top}
\end{align}
\begin{equation}\label{finite-xy}
[k_{11}^j {\bf U}_j^{x}]^y=\dfrac{\partial}{\partial y}\left ( k_{11}(x_i,y_j)\dfrac{\partial U(x_i,y_j,t)}{\partial x}\right ) 
\approx
D ({\bf K}^1_{j}{{{\bf U}}^y_j})
\doteq {\bf A}^{1}_j{\bf U}_j + {\bf b}^{x}_j, \,j=0,\dots, n.
\end{equation}
\begin{align}
{\bf A}^{1}_j &= d_0\beta^0_je_1^T +D_1\check{\bf K}^1_{j}D_2 -d_n\beta^n_je_{n+1}^T \quad \text{and}\\
{\bf b}^x_j &=\frac{1}{\cos^2 \theta(t)} \{ -d_0f_1(t,x_i,0)+ d_nf_2(t,x_i,1) \},
\end{align}
\item And more
\begin{align}\notag
[k_{22}^j {\bf U}_j^{x}]^y &\approx
\sum_{j=0}^n d_jk_{22}(x_i,y_j)  \dfrac{\partial 
U(x_i,y_j,t)}{\partial x } \\\notag
&= d_0k_{22}(0,y_j)  \dfrac{\partial U(x_0,y_j,t)}{\partial x } + 
\displaystyle \sum_{j=1}^{n-1} d_jk_{22}(x_i,y_j) 
\dfrac{\partial U(x_i,y_j,t)}{\partial x } 
 +d_nk_{22}(1,y_j)  \dfrac{\partial U(x_n,y_j,t)}{\partial x} \\\notag 
 &= \left [ d_0\alpha^0_j e_1^T + D_1\check{\bf K}^2_jD_2 -d_n\alpha^n_je_{n+1}^T\right ] 
{\bf U}_j+\frac{1}{\sin^2 \theta(t)} \{ -d_0 f_1(x_i,y_j,t) 
 + d_nf_2(x_i,y_j,t)\}\\ \label{bott-top}
\end{align}
\begin{equation}\label{finite-xy}
[k_{22}^j {\bf U}_j^{x}]^y=\dfrac{\partial}{\partial y}\left ( k_{22}(x_i,y_j)\dfrac{\partial U(x_i,y_j,t)}{\partial x}\right ) 
\approx
D ({\bf K}^2_{j}{{{\bf U}}^y_j})
\doteq {\bf A}^{2}_j{\bf U}_j + {\bf B}^{y}_j, \,j=0,\dots, n.
\end{equation}
\begin{align}
{\bf A}^{2}_j &= d_0\alpha^0_je_1^T +D_1\check{\bf K}^2_{j}D_2 -d_n\alpha^n_je_{n+1}^T \quad \text{and}\\
{\bf B}^x_j &=\frac{1}{\sin^2 \theta(t)} \{ -d_0f_1(t,x_i,0)+ d_nf_2(t,x_i,1) \},
\end{align}
\begin{align}\notag
{ \alpha}^0_j= 
\frac{1}{\sin^2 \theta(t)}  \left\{1-\frac{1}{2}[k_{11}(0,y_j)-k_{22}(0,y_j)] \sin(2 \theta(t))  -k_{11} (0,y_j)\cos^2 (\theta (t) \right\}
\end{align}
\begin{align}\notag
{\bf \alpha}^n_j= 
\frac{1}{\sin^2 \theta(t)}  \left\{1 -\frac{1}{2}[k_{11}(1,y_j)-k_{22}(1,y_j)]  \sin(2 \theta(t)) -k_{11} (1,y_j)\cos^2 (\theta (t)) \right\}
\end{align}

\begin{align}
{\bf A}^2_j &= \left [d_0 {\alpha}^0_je_1^T +D_1\check{\bf K}_j^{1}D_2 -d_n {\alpha}^n_j e_{n+1}^T\right ]^T 
\end{align}
\item 
\begin{align}\notag
[k_{11}^i {\bf U}_i^{y}]^x &\approx
\sum_{j=0}^n d_jk_{11}(x_i,y_j) U_y(t,x_i,y_j) \\\notag
&= d_0k_{11}(x_i,0) U_y(t,x_i,0) + 
\displaystyle \sum_{j=1}^{n-1} d_jk_{11}(x_i,y_j) 
 U_y(t,x_i,y_j)  
 +d_nk_{11}(x_i,1)  U_x(t,x_i,1)\\ 
&= \left [ d_0 \beta_i^0e_1^T +D_1\check{\bf K}_i^{1}D_2 -d_n \beta_i^ne_{n+1}^T\right ]^T
{\bf U}_i + {\bf b}^x_i\\ \label{bott-top}
&\doteq {\bf A}^1_i{\bf U}_i + {\bf b}_i^x,
\end{align}
$${\bf b}_i^x =\frac{1}{\sin^2 \theta(t)} \{ -d_0f_3(t,x_1,0)+ d_nf_4(t,x_1,1)\},$$
\begin{align}\notag
{ \beta}^0_i= 
\frac{1}{\sin^2 \theta(t)}  \left\{1-\frac{1}{2}[k_{11}(x_i,0)-k_{22}(x_i,0)] \sin(2 \theta(t))  -k_{22} (x_i,0)\cos^2 (\theta (t) \right\}
\end{align}
\begin{align}\notag
{\bf \beta}^n_i= 
\frac{1}{\sin^2 \theta(t)}  \left\{1 -\frac{1}{2}[k_{11}(x_i,1)-k_{22}(x_i,1)]  \sin(2 \theta(t)) -k_{22} (x_i,1)\cos^2 (\theta (t)) \right\}
\end{align}

\begin{align}
{\bf A}^1_i &= \left [d_0 {\beta}^0_ie_1^T +D_1\check{\bf K}_i^{1}D_2 -d_n {\beta}^n_i e_{n+1}^T\right ]^T 
\end{align}

\end{enumerate}

For a more detailed analysis on the CPM itself and theoretical aspects, the reader is again redirected to \cite{Bo_Luchesi_Bazan, B_Luchesi_Bazan}.

%
The discretization along all grid points with $U(x_i,y_j,t)$ enumerated in lexocographic order leads to  ${\bf U}$ corresponds to the block-wise vector
\begin{equation}
{\bf U} = \left [ \begin{array}{c}
{\bf U}_0 \\
{\bf U}_1 \\
\vdots \\
{\bf U}_n\\ \end{array} \right ],
\end{equation}
with every block defined in (\ref{Tj}) and ${\bf U}^0$ is the evaluation of initial condition $U_0$ from equation (\ref{eqe:Initial}) following the same ordering as ${\bf U}$. To organize and rewrite the eq. \ref{eq:heatF} in matricial form, we will create the following matrix
\begin{equation}\label{matrix_M}
{\bf M} =  \left ( \begin{array}{cccc}
{\bf F}_0 & & & \\
&  {\bf F}_1 & &  \\
& & \ddots &  \\
&   & & {\bf F}_{n}
\end{array} \right ) +    {\cal P} \left (\begin{array}{cccc}
{\bf G}_0 & & & \\
&  {\bf G}_1 & &  \\
& & \ddots &  \\
&   & & {\bf G}_{n}   
\end{array} \right )  {\cal P}^T ,
\end{equation}

where ${\bf F_j}(t)$ and ${\bf G_i}(t)$ are $1 \times n$ matrix given by
$${\bf F_j}=C^2[{\bf A}^1_j + {\bf b}^x_j]+Si^2[{\bf A}^{2}_j{\bf U}_j + {\bf b}^{y}_j]+Si[{\bf A}^{1}_j + {\bf b}^{x}_j]-Si[{\bf A}^{2}_j{\bf U}_j + {\bf b}^{y}_j]$$
$$=(C^2+Si){\bf A}^1_j+(Si^2-Si){\bf A}^{2}_j+(S_i+C^2){\bf b}^x_j+(Si^2-Si){\bf b}^{y}_j$$
Because, by eq. \ref{eq:heatF} we have
 $${\bf F_j}U_j(t)=C^2[k_{11}^j {\bf U}_j^{x}]^x(t)+Si^2[k_{22}^j {\bf U}_j^{x}]^x(t)+Si[k_{11}^j {\bf U}_j^{x}]^y(t)-Si[k_{22}^j {\bf U}_j^{x}]^y(t)$$
For Partial derivative in y
$${\bf G_i}=-Si[{\bf A}^{2}_i + {\bf b}^{y}_i] +Si[{\bf A}^1_i{\bf U}_i + {\bf b}_i^x]+Si^2[{\bf A}^1_i{\bf U}_i + {\bf b}_i^x]]+C^2[{\bf A}^2_j + {\bf b}^x_j]$$
$$=(C^2-Si){\bf A}^{2}_i+(Si+Si^2){\bf A}^1_i +(Si^2 +C^2){\bf b}^y_i+ (Si+Si^2) {\bf b}_i^x$$

By eq. \ref{eq:heatF} we have
$${\bf G_i}U_i(t)=-Si[k_{22}^i {\bf U}_i^{y}]^x(t)+Si[k_{11}^i {\bf U}_i^{y}]^x(t) +Si^2[k_{11}^i {\bf U}_i^{y}]^y(t)+C^2[k_{22}^i {\bf U}_i^{y}]^y(t)$$
Where
$$ C^2(t)_{1 \times n}=[ \cos^2 \theta(t)]$$
$$Si(t)_{1 \times n}=\frac{1}{2}[ \sin 2\theta(t)]$$
$$ Si^2(t)_{1 \times n}=[ \sin^2 \theta(t)]$$



\section{Time dependent System}
Therefore, discretization of equation \ref{eq:heatF} generate the following matrix system differential equation time dependent
\begin{equation}\label{disc_model}
\left \{ \begin{array}[l]{cl} 
  \textbf{U}'(t)=& {\bf H}(t,{\bf U})\nonumber\\
 \textbf{U}(0)=& \textbf{U}_0  
\end{array} \right.\end{equation}
\begin{equation}\label{model_disc}
{\bf H}(t,{\bf U}) =   {\bf M}{\bf U}(t)  + {\bf  Sg}(t),
\end{equation}
\begin{equation*}
{\bf Sg}(t) = S^x(t) + {\cal P}S^y(t),
\end{equation*}

for which we emphasize $S^x(t)$ and $S^y(t)$ have block entries ${\bf S}_j^x$ and ${\bf S}_j^y$ respectively, and $S(t)$
accommodates values   of the source term $g(t,x,y)$ along the grid, i.e, the block-entries of $S(t)$ are
\begin{equation}\label{disc_rhs}
S_j(t) = [g(t,x_0,y_j),\dots, g(t,x_n,y_j)]^T,\; j=0,1,\dots, n.
\end{equation}

Next step is to apply a numerical method for IVPs to (\ref{disc_model})  such as Euler explicit, Euler implicit, trapezoidal rule (Crank-Nicolson), Runge-Kutta class methods, or any other. We consider Crank-Nicolson's method (CN) 
for solving the IVP (\ref{disc_model}). Let a natural $M>0$ and a time discretization in the form $t_i = i\Delta t$, for $i = 0,1,\dots,M$, where $\Delta t = t_f/M$ is the time step, so that we have a uniformly spaced mesh on $[0,t_f]$. It's clear, though, that one can compute only a few iterations, i.e, taking $i$ from 0 to a natural smaller than $M$, generating solutions until some specific point of interest in time, as we do during numerical applications in the next sections. So, CN generates an approximate solution at time $t_{i+1}$ given by
 \begin{equation*}
 {\bf U}^{(i+1)} =  {\bf U}^{(i)} + \frac{\Delta t}{2} \left( {\bf M} {\bf U}^{(i)} +  
 {\bf  Sg}(t_{i}) + {\bf M} {\bf U}^{(i+1)} +  {\bf  Sg}(t_{i+1}) \right), \;\; i = 0,1,2,\dots, 
 \end{equation*}

\section{Inverse problem development}\label{sec:inverse}

Consider the vector of unknowns:
\begin{equation}\label{k_ortho}
{\bf k} = \left[ \begin{array}{c}
{\bf k}^{11}\\ 
{\bf k}^{22}
\end{array}  \right], \qquad  {\bf k}^{11}, \ {\bf k}^{22} \in \R^{(n+1)^2}, 
\end{equation}
where ${\bf k}^{11}$ and ${\bf k}^{22}$ contain a lexicographic enumeration of the wanted quantities $k_{11} (x_i,y_j)$ and $k_{22} (x_i,y_j)$ in Chebyshev's mesh. The objective is to recover conductivities $k_{11}$ and $k_{22}$ from the measurements of temperature at different time steps. Then, our problem consists basically on finding ${\bf k}^*$ such that
\begin{equation}\label{nonlin2}
{\bf k}^* = \argmin_{{\bf k} \in \R^{2(n+1)^2}} \phi({\bf k}), \quad \phi({\bf k}) =
\frac{1}{2} \|{\sf U}({\bf k}) - \widetilde{\sf U}\|_2^2,
\end{equation}
with $\widetilde{\sf U}$ representing the vector of known temperature. Moreover, ${\sf U}({\bf k})$ is the computed temperature relative to a given $\bf k$, i.e, we solve the direct problem for a certain $\bf k$ (which contains values for the conductivity only on mesh's points).  The nonlinear problem (\ref{nonlin2}) will be solved by Levenberg-Marquardt method (for nonlinear minimization) implemented by Yamashita and Fukushima 
and Morozov discrepancy principle
(for stopping criteria, due to the presence noise within data). Our version of LMM needs, at every step $j$, the computation of direction $d^j$ given by the equation
\begin{equation*}
(\mathcal{F}'^t({\bf k}^j) \mathcal{F}'({\bf k}^j) + \mu^j \textbf{D}^t \textbf{D}) 
{\bf d}^j = -\mathcal{F}'^t({\bf k}^j) \mathcal{F} ({\bf k}^j),
\end{equation*}
where $\bf D$ is an appropriate scaling matrix 
and $\mathcal{F}'({\bf k})$ denotes the derivative of $\sf T$ with respect to $\bf k$. Note that such derivative is not obviously computed, demanding an extra analysis of the model. In fact, in the continuous setting, variation of $\bf T$ with respect to $\bf k$ at a certain time step $t$ is given by the Jacobian
\begin{align*}
{\bf J}(t) &= \left [\frac{\partial {\bf T}}{\partial {\bf k}_1} ({\bf k},t), \dots, 
\frac{\partial {\bf T}}{\partial {\bf k}_{2(n+1)^2}}({\bf k},t)\right ]\\
& = \left [\frac{\partial {\bf T}}{\partial {\bf k}^{11}_1}({\bf k},t), \dots, 
\frac{\partial {\bf T}}{\partial {\bf k}^{11}_{(n+1)^2}} ({\bf k},t),
\frac{\partial {\bf T}}{\partial {\bf k}^{22}_1}({\bf k},t), \dots, 
\frac{\partial {\bf T}}{\partial {\bf k}^{22}_{(n+1)^2}}({\bf k},t)\right ].
\end{align*}
Then, the Jacobian is computed by assuming continuity of $\bf T$ and taking derivatives of (\ref{disc_model}) with respect to ${\bf k}_\ell$, $\ell \in \{ 1,2,\dots, 2(n+1)^2 \}$, resulting on the initial value problem for the $\ell$-th column of ${\bf J}(t)$,
\begin{equation}\label{jac-system}
\left \{
\begin{array}{l}

 \dfrac{\partial }{\partial t} \left( \dfrac{\partial {{\bf U}}}{\partial {\bf k}_{\ell}}({\bf k},t)\right) = 
{\bf M} \dfrac{\partial {{\bf U}}}{\partial {\bf k}_{\ell}}({\bf k},t)   + {\bf W}_{\ell}(t)     \\

\dfrac{\partial {\bf U}}{\partial {\bf k}_{\ell}} ({\bf k},0) = 0,
\end{array} \right .
\end{equation}
where $\bf M$ is the matrix between brackets in (\ref{model_disc}) and
\begin{equation}\label{jac-source}
{\bf W}_{\ell} (t) = \frac{\partial {\bf M}}{\partial {\bf k}_{\ell}}{\bf T}({\bf k},t), \;\ \ell =1, 2, \dots, 2(n+1)^2,
\end{equation}
are known as source terms. To write down such terms, we start by completing $D_1$ with a null vector in every side defining the matrix 
$$D_1^0 = \left[ {\bf 0}, D_1, {\bf 0} \right] \in \R^{(n+1) \times (n+1)}$$
and consequently denoting
\begin{equation}\label{hat_D}
\widehat{D} = \left[ \begin{array}{c}
e_1^T D_1^0 \\
D_2\\
e_{n+1}^T D_1^0
\end{array} \right] =: \left[ \hat{d}_0, \hat{d}_1, \dots, \hat{d}_n \right], \quad \hat{d}_i \in \R^{n+1}, \ i=0,1,\dots,n.
\end{equation}
Observe $\hat{d}_i \neq d_i$, $i = 0,1,\dots,n$. Then, due to the discretization model, one can verify the block-form of the source terms as follows:
\begin{equation*}
\frac{\partial {\bf M}}{\partial {\bf k}^{11}_1} = \left ( \begin{array}{clll}
\hat{d}_0 r_0^t & & & \\
&  {\bf 0}& &  \\
& & \ddots &  \\
&   & & {\bf 0}   \\
\end{array} \right ), \qquad 
\frac{\partial {\bf M}}{\partial {\bf k}^{22}_1} = {\bf P} \left ( \begin{array}{clll}
\hat{d}_0 r_0^t & & & \\
&  {\bf 0}& &  \\
& & \ddots &  \\
&   & & {\bf 0}   \\
\end{array} \right ) {\bf P},
\end{equation*}

\begin{equation*}
\frac{\partial {\bf M}}{\partial {\bf k}^{11}_{2}} = \left ( \begin{array}{clll}
\hat{d}_1 r_1^t & & & \\
&  {\bf 0} & &  \\
& & \ddots &  \\
&   & & {\bf 0}   \\
\end{array} \right ), \qquad 
\frac{\partial {\bf M}}{\partial {\bf k}^{22}_2} = {\bf P} \left ( \begin{array}{clll}
{\bf 0} & & & \\
& \hat{d}_0 r_0^t & &  \\
& & \ddots &  \\
&   & & {\bf  0}    \\
\end{array} \right )  {\bf P}^t,
\end{equation*}

\begin{equation*}
\vdots
\end{equation*}

\begin{equation*}
\frac{\partial {\bf M}}{\partial {\bf k}^{11}_{n+1}} = \left ( \begin{array}{clll}
\hat{d}_{n} r_n^t & & & \\
&  {\bf  0}& &  \\
& & \ddots &  \\
&   & & {\bf  0}    \\
\end{array} \right ), \qquad \frac{\partial {\bf M}}
{\partial {\bf k}^{22}_{n+1}} = {\bf P} \left ( \begin{array}{clll}
{\bf  0} & & & \\
& \ddots & &  \\
& & {\bf  0} &  \\
&   & & \hat{d}_0 r_0^t    \\
\end{array} \right )  {\bf P}^t,
\end{equation*}

\begin{equation*}
\frac{\partial {\bf M}}{\partial {\bf k}^{11}_{n+2}} = \left ( \begin{array}{clll}
{\bf 0} & & & \\
&  \hat{d}_0 r_0^t& &  \\
& & \ddots &  \\
&   & & {\bf 0}   \\
\end{array} \right ), \qquad 
\frac{\partial {\bf M}}{\partial {\bf k}^{22}_{n+2}} = {\bf P} \left ( \begin{array}{clll}
\hat{d}_1 r_1^t & & & \\
&  {\bf 0}& &  \\
& & \ddots &  \\
&   & & {\bf 0}   \\
\end{array} \right ) {\bf P},
\end{equation*}

\begin{equation*}
\frac{\partial {\bf M}}{\partial {\bf k}^{11}_{n+3}} = \left ( \begin{array}{clll}
{\bf 0} & & & \\
&  \hat{d}_1 r_1^t & &  \\
& & \ddots &  \\
&   & & {\bf 0}   \\
\end{array} \right ), \qquad 
\frac{\partial {\bf M}}{\partial {\bf k}^{22}_{n+3}} = {\bf P} \left ( \begin{array}{clll}
{\bf 0} & & & \\
& \hat{d}_1 r_1^t & &  \\
& & \ddots &  \\
&   & & {\bf  0}    \\
\end{array} \right )  {\bf P}^t,
\end{equation*}

\begin{equation*}
\vdots
\end{equation*}

\begin{equation*}
\frac{\partial {\bf M}}{\partial {\bf k}^{11}_{2(n+1)}} = \left ( \begin{array}{clll}
{\bf  0} & & & \\
&  \hat{d}_{n} r_n^t & &  \\
& & \ddots &  \\
&   & & {\bf  0}    \\
\end{array} \right ), \qquad \frac{\partial {\bf M}}
{\partial {\bf k}^{22}_{2(n+1)}} = {\bf P} \left ( \begin{array}{clll}
{\bf  0} & & & \\
& \ddots & &  \\
& & {\bf  0} &  \\
&   & & \hat{d}_1 r_1^t    \\
\end{array} \right )  {\bf P}^t,
\end{equation*}
and forth, replicating the rule until index $(n+1)^2$. Observe that to compute the full Jacobian ${\bf J} (t)$ it's necessary solving $2(n+1)^2$ IVPs, since each column is represented by a system such as (\ref{jac-system}). Finally, the sought derivative is a collection of Jacobians at different time steps, namely
\begin{equation}\label{jacobian}
\mathcal{F}'({\bf k}) = \left[ \begin{array}{c}
{\bf J}_1\\ 
{\bf J}_2\\ 
\vdots \\ 
{\bf J}_N
\end{array}  \right], 
\end{equation}
where ${\bf J}_q = {\bf J}(t_q)$ represents the Jacobian of $\bf U$ with respect to $\bf k$ 
at time step $t_q$, for $q = 1,2,\dots, N$, where $N$ represents the number of time steps under consideration.

Then, the inverse problem is attacked with the known strategy: minimization of residual temperature through Levenberg-Marquardt method along with a stopping rule (discrepancy principle, relative residual, etc.) to control noise dispersion. In next subsection we show some of the obtained results as well as comments on the stopping criterion.

\section{Conclusion}
The reconstruction process involves solving the governing partial differential equation (PDE) using a semi-discrete approach, which integrates a pseudospectral collocation method for spatial discretization and the numerical scheme for temporal discretization. Additionally, to address the numerical inversion problem, the two-dimensional transient diffusion  is reconstructed through a nonlinear least squares optimization, employing the Levenberg-Marquardt method (LMM), with the Morozov discrepancy principle as a stopping criterion to mitigate the effects of noise in the data.
To demonstrate the efficacy of this pseudospectral approach, the CPM will be specifically applied to the reconstruction of the diffusion tensor in some benchmark tests and generate synthetic experimental data. 
  \bibliographystyle{plain}

\end{document}